\documentclass[12pt,reqno]{amsart} 
\usepackage{amsthm,amsfonts,amssymb,amscd}

\newtheorem{theorem}[subsection]{Theorem}
\newtheorem{proposition}[subsection]{Proposition}

\newtheorem{lemma}[subsection]{Lemma}
\newtheorem{corollary}[subsection]{Corollary}

\theoremstyle{definition}

\newtheorem{proposition-definition}[subsection]{Proposition-Definition}

\theoremstyle{remark}
\newtheorem{remark}[subsection]{Remark}

\newcommand{\Hff}{H_{14,15}}
\newcommand{\Hfo}{\makebox[1em]{\raisebox{8pt}[2pt]{$\scriptstyle\circ$}}\hspace{-1em}{H}_{14,15}}
\newcommand{\Mfo}{\makebox[1.1em]{\raisebox{6.6pt}[1pt]{$\scriptstyle\circ$}}\hspace{-1.3em}{\mathcal M}^4_{15}}
\newcommand{\dual}{{\scriptscriptstyle\vee}}
\newcommand{\mt}[1]{\operatorname{#1}}

\newcommand{\codim}{\operatorname{codim}\nolimits}
\newcommand{\im}{\operatorname{im}\nolimits}
\newcommand{\rk}{\operatorname{rk}\nolimits}

\newcommand{\corank}{\operatorname{corank}\nolimits}

\newcommand{\Ext}{\operatorname{Ext}\nolimits}

\newcommand{\Aut}{\operatorname{Aut}\nolimits}

\newcommand{\Cl}{\operatorname{Cl}}
\newcommand{\Pf}{\operatorname{Pf}\nolimits}

\newcommand{\Pic}{\operatorname{Pic}\nolimits}

\newcommand{\End}{{{\mathcal E}nd\:}}

\newcommand{\PD}{\partial}

\newcommand{\CC}{{\mathbb C}}

\newcommand{\ZZ}{{\mathbb Z}}

\newcommand{\PP}{{\mathbb P}}

\newcommand{\OOO}{{\mathcal O}}

\newcommand{\III}{{\mathcal I}}

\newcommand{\EEE}{{\mathcal E}}

\newcommand{\LLL}{{\mathcal L}}
\newcommand{\FFF}{{\mathcal F}}
\newcommand{\GGG}{{\mathcal G}}

\newcommand{\KKK}{{\mathcal K}}
\newcommand{\NNN}{{\mathcal N}}
\newcommand{\MMM}{{\mathcal M}}

\newcommand{\mapto}[1]{\stackrel{#1}{\longrightarrow}}

\newcommand{\isoto}{\stackrel{\textstyle\sim}{\longrightarrow}}
\newcommand{\map}[1]{\stackrel{#1}{\longrightarrow}}

\newcommand\alp{\alpha}

\newcommand\eps{\epsilon}

\newcommand\si{\sigma}

\newcommand\lra{{\longrightarrow}}
\newcommand\rar{\rightarrow}

\author{A. Iliev}

\address{A. I.: Inst. of Math.,
Bulgarian Acad. of Sci.,
Acad. G. Bonchev Str., 8,\ 
1113 Sofia, Bulgaria}
\email{ailiev@math.bas.bg}

\author{D. Markushevich}

\address{D. M.: Math\'ematiques - b\^{a}t. M2, Universit\'e Lille 1,
F-59655 Villeneuve d'Ascq Cedex, France}
\email{markushe@gat.univ-lille1.fr}


\title{Quartic 3-fold: Pfaffians, instantons
and half-canonical curves}

\begin{document}

\begin{abstract} A generic quartic 3-fold $X$ admits a 7-dimensional family
of representations as the Pfaffian of an 8 by 8 skew-symmetric matrix
of linear forms. This provides a 7-dimensional
moduli space $M$ of rank 2 vector bundles on $X$.
A precise geometric description of a 14-dimensional family of
half-canonical curves $C$ of genus 15 in $X$ such that
the above vector bundles are obtained by Serre's construction from $C$
is given. It is proved that the Abel--Jacobi map of this family factors through
$M$, and the resulting map from $M$ to the intermediate Jacobian is
quasi-finite. In particular, every component of $M$ has non-negative
Kodaira dimension. Some other constructions of rank 2 vector bundles with
small Chern classes are discussed; it is proved that the smallest
possible charge of an instanton on $X$ is 4.
\end{abstract}

\maketitle

\section*{Introduction}

This paper is a part of
the study of moduli spaces of vector bundles
with small Chern classes on certain Fano threefolds.
It provides some non-existence results and constructions of a few moduli components of vector bundles
on a quartic threefold. One of them is the component of kernel bundles,
defined similarly to
that of \cite{MT}, \cite{IM} for the case of a cubic threefold.
Our work received a strong pulse with the publication of the paper
of Beauville \cite{B},
which allowed to simplify some arguments used in \cite{MT} and at the
same time put our results in a more general framework
of Pfaffian hypersurfaces. We also prove that there are no normalized rank 2 stable
vector bundles on $X$ with $c_2<4$ and exhibit two constructions of stable
vector bundles 
with $c_1=0$ and $c_2=4$. Only one of them provides instantons; 
we show in fact
that all the instantons of charge 4 are obtained by this construction.

In \cite{MT}, it was proved that the 
Abel--Jacobi map of the family of normal elliptic quintics lying
on a general cubic threefold $V$ factors through a moduli component of
stable rank 2 vector bundles on $V$
with Chern numbers $c_1=0,c_2=2$,
whose general point represents
a vector bundle obtained by Serre's construction from an elliptic
quintic. The elliptic quintics mapped to a point of the moduli
component vary in a 5-dimensional projective space
inside the Hilbert scheme of curves,
and the map from the moduli component to the intermediate Jacobian
is quasi-finite. Later, in \cite{IM}, this modular component was identified
with the variety of representations of $V$ as a linear section
of the Pfaffian cubic in $\PP^{14}$ and it was proved that the
degree of the quasi-finite map is 1, so the moduli component is birational
to the intermediate Jacobian $J^2(X)$. Beauville mentions in \cite{B}
a recent work of Druel, yet unpublished, which proves
that the moduli space $M_V(2;0,2)$ is irreducible, so
its unique component is the one described above.

In the present paper, we prove that a generic quartic threefold
$X$ admits a 7-dimensional family of essentially different
representations as the Pfaffian of an $8\times 8$ skew-symmetric matrix
of linear forms. Thanks to \cite{B}, this provides a 7-dimensional
family of arithmetically Cohen--Macaulay (ACM for short) vector bundles
on $X$, obtained as the bundles of kernels of the
$8\times 8$ skew-symmetric matrices of rank 6 representing points of $X$.
We show that this family is a smooth open set $M_X$ in the
moduli space of stable vector bundles $M_X(2;3,14)\simeq M_X(2;-1,6)$.
The ACM property means the vanishing of the
intermediate cohomology $H^i(X,\EEE (j))$ for all $i=1,2$, $j\in\ZZ$.

We give also a precise geometric characterization of the ACM
curves arising as schemes of zeros of sections of the above kernel
vector bundles. According to Beauville, they are half-canonical ACM
curves of degree 14 in $\PP^4$; we show that they are linear
sections of the rank 4 locus $Z\subset\PP (\wedge^2\CC^7)$ in the
projectivized space of the $7\times 7$ skew-symmetric matrices.
Linear sections of $Z$ arose already in the literature: R\o dland
\cite{R} studied the sections $\PP^6\cap Z$, which are
Calabi--Yau threefolds. We show that such curves fill out open sets
of smooth points of the Hilbert schemes of $X$ (of dimension 14)
and of $\PP^4$ (of dimension 56), and that the isomorphism
classes of smooth members of this family fill out a 32-dimensional
moduli component $\Mfo$ of curves of genus 15 with a 
theta-characteristic linear series of dimension $4$.

Next we study the Abel--Jacobi map of the ACM half-canonical curves
of genus 15 in $X$. It factors through $M_X$ via Serre's construction:
the fibers over points of $M_X$ are $\PP^7$, and the resulting map
from $M_X$ to $J^2(X)$ is \'etale quasi-finite, hence its image is 7-dimensional.
The role of the above half-canonical curves is similar to that of
normal elliptic quintics in the case of the cubic threefold $V$, where
the Abel--Jacobi map factors through the instanton moduli space
with fibers $\PP^5$ and with a 5-dimensional image; as $\dim J^2(V)=5$, 
the image is an open subset of $J^2(V)$.
The result for a quartic threefold is somewhat weaker: here we do not
know whether the degree of the quasi-finite map is 1 and whether $M_X$
is irreducible.
Moreover, as $7=\dim M_X<30=\dim J^2(X)$, we cannot conclude,
as in the case of a cubic threefold, that the image of $M_X$ is
an open subset of an Abelian variety; we can only state
that every component of it, and hence of $M_X$ itself, has a non-negative
Kodaira dimension.

In Section 1, we gather preliminary results: a criterion for stability
of rank 2 sheaves, Bogomolov's inequality, and prove the non-existence
of normalized rank 2 stable vector bundles with small second Chern classes.
We provide two constructions of such vector bundles for $c_1=0$, $c_2=4$
and discuss the geometry of the instanton component(s).

In Section 2, we prove that a generic quartic 3-fold is Pfaffian,
in using the same method as was used by Adler in his Appendix to
\cite{AR} for a cubic threefold: take a particular quartic
which is Pfaffian and prove that the differential of the
Pfaffian map from the family of all the $8\times 8$ skew-symmetric
matrices of linear forms to the family of quinary quartics
is of maximal rank.
We prove also basic facts about $M_X$:
stability, dimension 7, smoothness.

Section 3 treats half-canonical ACM curves 
of genus 15 on $X$ and in~$\PP^4$.

Section 4 applies the technique of the Tangent Bundle Sequence
of Clemens--Griffiths \cite{CG} and Welters \cite{W} to the
calculation of the differential of the Abel--Jacobi map
for the family of the above half-canonical curves $C$.
It identifies the kernel of the differential with
$H^1(\NNN_{C/\PP^4}(-1))^\dual$, and we prove that it has dimension 7.

\smallskip

{\em Acknowledgements}. The second author acknowledges with pleasure the hospitality
of the Max-Planck-Institut f\"ur Mathematik at Bonn, where
he wrote this paper.

\section{Generalities and the case \( c_1=0 \)}

Let $X$ be a smooth quartic threefold. It is well known that $\Pic (X)$ is
isomorphic to $\ZZ$, generated by the class of the hyperplane section $H$, and
the group of algebraic $1$-cycles modulo topological equivalence is also
isomorphic to $\ZZ$, generated by the class of a line $l\subset X$.
For two integers $k,\alp$, we will denote by $M_X(2;k,\alp )$ the moduli space
of {\em stable} rank 2 vector bundles $\EEE$ on $X$ with Chern classes $c_1=k[H]$
and $c_2=\alp [l]$. We will often identify the Chern classes with integers in
using the generators $[H], [l]$ of the corresponding groups of algebraic cycles.
We have $[H]^2=4[l]$.

By the definition of the Chern classes \sloppy and by Riemann--Roch-- Hirzebruch,
we have for $\EEE\in M_X(2;k,\alp )$:

$
c_1(\EEE (n))=c_1(\EEE )+2n[H]=(k+2n)[H],\ c_2(\EEE (n))=c_2(\EEE )+n[H]c_1(\EEE )
+n^2[H]^2=(\alp +4kn+4n^2)[l],\ 
\chi (\EEE )=\frac{2}{3}k^3-\frac{1}{2}k\alp +k^2-\frac{1}{2}\alp
+\frac{7}{3}k +2\ .
$

A rank 2 torsion free sheaf $\EEE$ on $X$ is {\em normalized}, if $c_1(\EEE )=
k[H]$ with $k=0$ or $k=-1$.
We can make $\EEE$ normalized in replacing it by a suitable twist $\EEE (n)$.
If $\EEE$ is locally free, we have $\EEE^\dual\cong\EEE\otimes (\det\EEE )^{-1}$,
so that $\EEE$ is self-dual when $k=0$.

The following lemmas are well known:

\begin{lemma}\label{criterion}
Let $\EEE$ be a normalized rank 2 reflexive sheaf on a nonsingular projective
variety $X$ with $\Pic (X)\simeq\ZZ$. Then it is stable
if and only if $h^0(\EEE )=0$.
\end{lemma}

\begin{proof}
Any saturated torsion free rank 1
subsheaf of $\EEE$ 
is invertible of the form $\OOO_X(m)$
and gives an exact triple
\begin{equation}\label{diese}
%$$
0\lra \OOO_X(m)\lra \EEE \lra \III_{Z}(k-m)
\lra 0\ ,
%\eqno(16)
%$$
\end{equation}
where $Z$ is a subscheme of $X$ of codimension 2.
The triple breaks the Gieseker stability
of $\EEE$ if and only if $m\geq 0$.

If we assume that $\EEE$ has global sections, then there exists a
triple (\ref{diese}) with $m=0$, hence $\EEE$ is not stable.
If we assume that $\EEE$ has no global sections, then in any triple
(\ref{diese}) for $\EEE$, we have $m< 0$, because $h^0(\EEE )
\geq h^0(\OOO_X(m))$. 
Hence $\EEE$ is stable.
\end{proof}

\begin{lemma}\label{bogomolov}
Let $\FFF$ be a rank $r$ semistable reflexive sheaf on $X$. Then
$(2rc_2(\FFF )-(r-1)c_1^2(\EEE ))\cdot H\geq 0$. If $r=2$ and
$\FFF$ is stable, then the inequality is strict.
\end{lemma}

\begin{proof}
This is Bogomolov's Theorem, proved by him for $T$-(semi)stable sheaves \cite{Bo}.
For another approach to the proof and for relations between different
notions of (semi)stability, see e. g. \cite{Ko}.
\end{proof}

\begin{proposition}
Let $X$ be a smooth quartic threefold. Then the following
statements hold:

(i) $M_X(2;0 ,\alp )=\varnothing$
for all odd $\alp$ and for $\alp\leq 2$.

(ii) $M_X(2;-1 ,\alp )=\varnothing$ for
$\alp \leq 3$.
\end{proposition}

\begin{proof}
The case of odd $\alp$ in (i) follows trivially from Riemann--Roch--Hirzebruch:
we have $\chi (\EEE )=2-\frac{1}{2}\alp$. For the remaining cases,
the proof goes exactly as in \cite{B-MR}. The first step is to show that if
$M_X(2;\eps ,\alp)\neq\varnothing$ ($\eps =0,-1$), then
$h^0(\EEE (1))\neq 0$ for all $\EEE\in M_X(2;\eps ,\alp)$. The second step
is to verify that there are no curves on $X$ that might be zero loci of
sections of $\EEE (1)$.

So, let $\EEE\in M_X(2;\eps ,\alp)$, $\eps =0,-1$, $\alp \leq 2-\eps$.
Assume that $h^0(\EEE (1))= 0$. By Serre duality, $h^3(\EEE (1))=h^0(\EEE (\eps -2))=0$.
Hence $h^2(\EEE (1))\geq\chi (\EEE (1))$. We have $\chi (\EEE (1))=10-\frac{3}{2}\alp$
if $\eps =0$ and $\chi (\EEE (1))=6-\alp$ if $\eps =-1$. Hence
$\dim\Ext^1(\EEE ,\OOO (-2))=h^1(\EEE^\dual (-2))> 0$, and there exists a
non-trivial extension of vector bundles
\begin{equation}\label{extension}
0\lra \OOO_X(-2)\map{\sigma} \FFF \lra \EEE\lra 0\ .
\end{equation}
We have $\Delta (\FFF)= (c_1(\FFF )^2-3c_2(\FFF ))H=16-3\alp >0$ if
$\eps =0$ and $\Delta (\FFF)=12-3\alp$ if $\eps =-1$, so $\FFF$ is unstable by
Lemma \ref{bogomolov}. To fix ideas, restrict ourselves to the case $\eps = 0$,
the other case being completely similar.

The unstability of $\FFF$ can manifest itself in two ways: either $\FFF$
contains a rank 1 saturated subsheaf $\OOO_X (n)$ with $n> c_1(\FFF )/\rk \FFF=-2/3$,
or there exists a non-trivial morphism of sheaves $\FFF\map{\phi}\OOO_X (n)$ with
$n< c_1(\FFF )/\rk \FFF$. In the first case, $n\geq 0$, hence $h^0(\FFF )\neq 0$,
hence $h^0(\EEE )\neq 0$, and this contradicts the stability of $\EEE$.

In the second case, $n\leq -1$. If $n<-2$, then $\phi\sigma =0$ and $\phi$
descends to a non-trivial morphism $\EEE\lra\OOO_X(n)$, which contradicts
the stability of $\EEE$. Hence $n=-1$ or $-2$. It cannot be $-2$, because
otherwise the extension (\ref{extension}) would be split. So $n=-1$ and we
obtain the exact triple
$$
0\lra \FFF '\lra\FFF\map{\phi}\III_Z(-1)\lra 0\ ,
$$
in which $\FFF '$ is a reflexive rank 2 sheaf. We have 
$c_1(\FFF ')=-[H]$, $c_2(\FFF ')\leq
c_2(\FFF )-c_1(\FFF ')c_1(\OOO_X(1))=(\alp -4)[l]<0$, hence $\FFF '$
is unstable by Bogomolov's inequality. By Lemma \ref{criterion},
$h^0 (\FFF ')\neq 0$, which implies $h^0 (\FFF )\neq 0$ and hence
$h^0 (\EEE )\neq 0$. This contradicts the stability of $\EEE$.

Thus we have proved $h^0 (\EEE (1))\neq 0$. As $\Pic X\cong\ZZ$,
the scheme $C=C_s$ of zeros of a non-trivial section $s$ of $\EEE (1)$
is a l. c. i. of pure codimension 2. Hence $\EEE (1)$ fits into the
following exact triple
\begin{equation}\label{serre-general}
0\lra\OOO_X\map{s}\EEE (1)\lra\III_C(2)\lra 0\ .
\end{equation}
We have $[C]=(4+\alp )[l]$, $\omega_C=\OOO_C(1)$, $2p_a(C)-2=4+\alp$.
It remains only to verify
the case of $\alp =2$. We have $p_a(C)=4$ and $C$ is embedded into $\PP^4$
by a subsystem of the canonical system. The exact triple (\ref{serre-general}),
twisted by $\OOO_X(-1)$, implies that $C$ is not contained in a
hyperplane. Hence $C$ is not connected. It cannot be a union of more than
one connected components either, because at least one of them should be
of degree $\leq 3$ and hence $\omega_C$ cannot be ample.

The proof is completed in a similar way in the case $k=-1$.
\end{proof}

For $k=0$, we have proved that there are no stable rank 2 vector
bundles with $\alp <4$. However, they do exist for $\alp =4$.
Indeed, assuming that $h^0(\EEE (1))\neq 0$, we find only two
possibilities for the zero locus $C$ af a non-trivial section of
$\EEE (1)$: either $C$ is a canonical curve of degree 8, or the union
of two canonical curves of degree 4 which spans $\PP^4$ (including the 
non-reduced limits of such curves). In the first case, standard
calculations show that the generic quartic
contains an 8-dimensional family of such curves $C$.
Only in this case the vector
bundle $\EEE$ has {\em natural cohomology}, that is, for every $t$,
$h^i(\EEE (t))\neq 0$ for at most one value of $i$. It is reasonable to call
{\em  instantons} the stable vector bundles with natural cohomology,
such that $c_1(\EEE )=0$ and the instanton
condition $h^1(\EEE(-2))=0$ is verified. Thus we have the following statement:

\begin{proposition}\label{Moduli}
(i) On any smooth quartic threefold $X$, there is an irreducuble
component $M_X^0$ of $M_X(2;0,4)$ which parametrizes the vector
bundles obtained by Serre's construction from the curves $C=C_1^4\sqcup C_2^4$,
where $C_i^4$ are plane sections of $X$. These vector
bundles satisfy $h^1(\EEE (1))=1$, hence they are not instantons.

(ii) Let $C$ be a smooth complete intersection of $3$ quadrics in $\PP^4$.
Then there exists a smooth quartic threefold $X$ containing $C$, and
the vector bundles on $X$ obtained
by Serre's construction (\ref{serre-general}) from the curve $C$ 
and from its generic deformations in $X$ sweep out a component
$M_X^1(C)$ of $M_X(2;0,4)$, different from $M_X^0$. The vector bundles
$\EEE\in M_X^1(C)$ are instantons.

(iii) Let $X$ be a smooth quartic threefold. Then any component $M$
of $M_X(2;0,4)$, such that $h^0(\EEE (1))\neq 0$ for some $\EEE\in M$,
is one of the above components $M_X^0$, $M_X^1(C)$.
\end{proposition}

\begin{remark}
Let ${\mathcal C}^5_8(X)$ be the 8-dimensional family of curves 
$C$ on the general $X = X_4$ as in Proposition \ref{Moduli}, (ii), 
and let $v_2:X \rightarrow \PP^{14}$ be the Veronese map. 
For $C \in {\mathcal C}^5_8(X)$, one can see that the \mbox{11-space} 
$\PP^{11}(C) =<\! v_2(C)\! > \subset \PP^{14}$ 
lies in a unique rank 6 quadric $Q = Q_C \supset v_2(X)$. Indeeed,
if $C$ is the intersection of three quadrics $q_i=0$ ($i=1,2,3$),
then the equation of $X$ can be written in the form
$q_1(x)\tilde{q}_1(x)+q_2(x)\tilde{q}_2(x)+q_3(x)\tilde{q}_3(x)=0$,
which provides the rank 6 quadric
$Q=\{ l_1\tilde{l}_1+l_2\tilde{l}_2+l_3\tilde{l}_3=0\}$,
where $l_i,\tilde{l}_i$ are the linear forms in the Veronese embedding
corresponding to the quadratic forms  $q_i,\tilde{q}_i$.
It is a degenerate cone whose ridge is $\PP^8$, the kernel
of the quadratic form, and whose base is $G$, 
a non-degenerate 4-dimensional
quadric in~$\PP^5$. The curves $C_s$ of zeros of sections
$s\in H^0(\EEE (1))$ form a projective space $\PP^3$, naturaly identified
with one of the two $\PP^3$'s parametrizing projective 2-planes 
$\{ l_1'=l_2'=l_3'=0\}$ in $G$, the equations of $C_s$ being of the
form $q_1'=q_2'=q_3'=0$, where $q_i'$ correspond to $l_i'$ under the Veronese
map. The above $\PP^{11}$ is the cone over this $\PP^3$ with ``vertex"~$\PP^8$.

Let ${\Lambda}_C = \{\PP^{11}_t: t \in {\PP}^3 \}$ be the 
ruling of $Q_C$ defined by ${\bf P}^{11}(C) \in \Lambda$, 
${\bf P}^{11}(C) = {\bf P}^{11}_o$. 
Then the sections of $E_C(1)$ are exactly the curves 
$C_t = {\bf P}^{11}_t \cap v_2(X) \subset v_2(X) \cong X$, 
$t \in {\bf P}^3$, $C = C_o$.

Thus we can represent $M_X^1$, the union of the 5-dimensional components
$M_X^1(C)$ of Proposition \ref{Moduli}, as the variety
$\tilde{D}_6(v_2(X)):=\{\:$the $\PP^{11}$'s contained in some rank 6 quadric $Q\supset v_2(X)\:\}$,
which is, in its turn, the double cover of the family 
$D_6(v_2(X)):=\{ \:$the rank 6 quadrics $Q\supset v_2(X)\:\}$.
\end{remark}

\begin{remark}
According to \cite{Tyu}, (3.1.45), the virtual dimension of $M_X(2;0,2)$
is 1, so one could expect a curve of isomorphism classes of stable vector bundles
with $c_1=0,c_2=2$. But we have proved that $M_X(2;0,2)$ is empty, providing thus
one more example of a situation when dimension is different from the virtual one. 
\end{remark}

\begin{remark}
For $k=-1$, we leave open the cases of $\alp =4,5$. We construct in what
follows a 7-dimensional component, or a union of 7-dimensional
components of $M_X(2;-1,6)$ (we do not approach the question on
the number of these components).
\end{remark}

\section{Generic quartic 3-fold is Pfaffian}

Let $E$ be an 8-dimensional vector space over $\CC$. Fix a basis
$e_0,\ldots ,e_7$
for  $E$, then $e_{ij}=e_i\wedge e_j$ for
$0\leq i<j\leq 7$ form a basis for the Pl\"ucker space $\wedge^2E$.  Let
$x_{ij}$ be the corresponding (Pl\"ucker) coordinates.
The embedding of the Grassmannian $G=G(2,E)$ in 
$\PP^{27}=\PP 
(\wedge^2E)$ is precisely the locus of rank 2 skew symmetric
$8\times 8$ matrices $M$ with elements $x_{ij}$ above the diagonal.
Let $G\subset\Omega\subset\Xi\subset \PP^{27}$ be the filtration
of $\PP^{27}$ by the rank of $M$, that is
$\Omega =\{ M\mid \rk M\leq 4\}$, $\Xi =\{ M\mid \rk M\leq 6\}$.
Then $G$, $\Omega\setminus G$, $\Xi\setminus\Omega$ and $\PP^{27}\setminus
\Xi$ are orbits of $PGL(8)$, acting via $\wedge^2$ of its
standard representation (see e. g. \cite{SK}),
and we have $G=\mt{Sing}\Omega$, $\dim G=12$,
$\Omega =\mt{Sing}\Xi$, $\dim\Omega = 21$. 
$\Xi$ is defined by the quartic equation 
$\Pf (M)=0$, where $\Pf$ stands for the Pfaffian of a skew-symmetric matrix.
We will call $\Xi$ the {\em Pfaffian hypersurface}
of~$\PP^{27}$.

Let $H\subset\PP^{27}$ be a 4-dimensional linear subspace which is
not contained in $\Xi$. Then the intersection $H\cap \Xi$
will be called a Pfaffian quartic 3-fold.
As $\codim_\Xi\Omega =5$, the linear section $H\cap \Xi$ is nonsingular
for general $H$. Suppose that a quartic 3-fold
$X\subset \PP^4$ has two different
representations $\phi_1:X\isoto H_1\cap\Xi$, $\phi_2:X\isoto H_2\cap\Xi$
as linear sections of $\Xi$. We will call them equivalent if 
$\phi_2\circ\phi_1^{-1}$ is the restriction of a transformation
from $PSL(8)$.

\begin{proposition}\label{a-la-Adler}
A generic quartic 3-fold admits a $7$-parameter family
of non-equivalent representations as linear sections
of the Pfaffian hypersurface in $\PP^{27}$.
\end{proposition}

\begin{proof}
 We are using the same argument as that of \cite{AR}, Theorem (47.3).
The family of quartic 3-folds in $\PP^4$ is parametrized by $\PP^{69}$,
and that of the Pfaffian representations of quartic 3-folds by an open set
in the variety $\mt{Lin}(\PP^4,\PP^{27})$ of linear morphisms between the
two projective spaces.
So, we are  going to specify one particular quartic 3-fold $X_0=\{ F_0=0\}$ which admits
a Pfaffian representation $F_0=\Pf (M_0)$, then we will show that the differential of the
map $f:\mt{Lin}(\PP^4,\PP^{27})\dasharrow \PP^{69}$ at $M_0$
is surjective, and this will imply that $f$ is dominant.

Let
$$ M_0=\left[\begin{array}{cccccccc} 0&x_1&x_2&x_3&x_4&x_5&x_1&0\\ 
     -x_1&0&0&x_5&0&0&-x_3&-x_1\\ 
     -x_2&0&0&x_1&x_1&0&0&-x_4\\ 
     -x_3&-x_5&-x_1&0&x_2&0&0&0\\ 
     -x_4&0&-x_1&-x_2&0&x_3&x_1&0\\ 
     -x_5&0&0&0&-x_3&0&x_4&x_2\\ 
     -x_1&x_3&0&0&-x_1&-x_4&0&x_5\\ 
     0&x_1&x_4&0&0&-x_2&-x_5&0
 \end{array}\right]
,$$

$
F_0=\Pf (M_0)=x_1^3x_2-x_1^3x_3+x_2^3x_3-x_1x_2x_3^2-x_1x_2^2x_4+
x_1^2x_3x_4+x_1x_2x_3x_4+x_3^3x_4-x_1^2x_4^2+x_1x_2x_4^2+x_1^3x_5-
x_1^2x_2x_5-x_1x_2^2x_5-x_1^2x_3x_5+x_1x_3x_4x_5+x_2x_3x_4x_5+x_4^3x_5+
x_2x_3x_5^2-x_1x_4x_5^2+x_1x_5^3. $

A point $M\in \mt{Lin}(\PP^4,\PP^{27})$ is the proportionality class of an $8\times 8$
skew-symmetric matrix of linear forms $l_{ij}$ and is given by its
$5\cdot 28=140$ homogeneous coordinates $(a_{ijk})$ such that
$l_{ij}=\sum_k a_{ijk}x_k$ ($0\leq i<j\leq 8,1\leq k\leq 5$).
We have $\PD f(M)/\PD a_{ijk}\ =\ x_k\Pf_{ij}(M)$,
where $\Pf_{ij}(M)$ denotes the Pfaffian of the $6\times 6$ matrix
obtained by deleting the $i$-th and the $j$-th rows and
the $i$-th and the $j$-th columns of $M$.

Computation by the Macaulay 2 program \cite{M2} shows that, for the
above matrix $M_0$, the 140 quartic forms $x_k\Pf_{ij}(M_0)$ generate
the whole 70-dimensional space of quinary quartic forms, hence
$f$ is of maximal rank at $M_0$. One can also easily make Macaulay 2
to verify that
$X_0$ is in fact nonsingular, though this is not essential for the
above proof.

It \sloppy remains to verify that the generic fiber of the induced
map $\bar{f}:{\raisebox{-2pt}{$PGL(5)$}\!\!\setminus\!\!
\raisebox{2pt}{$\mt{Lin}(\PP^4,\PP^{27})$}\! /
\!\raisebox{-2pt}{$PGL(8)$}}\dasharrow{\raisebox{-2pt}{$PGL(5)$}\!\!\setminus\!\!
\raisebox{2pt}{$\PP^{69}$}}$
is 7-dimensional. By counting dimensions, one sees that this is equivalent to
the fact that the stabilizer of a generic point of the Grassmannian
$G(5,28)={\raisebox{-2pt}{$PGL(5)$}\!\!\setminus\!\!
\raisebox{2pt}{$\mt{Lin}(\PP^4,\PP^{27})$}}$ in $PGL(8)$ is 0-dimensional.

Take a generic 4-dimensional linear subspace $H\subset\PP^{27}$. Then the
quartic 3-fold $X=H\cap\Xi$ is generic, and hence $\Aut (X)$ is trivial.
Thus the stabilizer $G_H$ of $H$ in $PGL(8)$ acts trivially an $X$, and
hence on $H$. This implies the triviality of $G_H$ by (5.3) of \cite{B}.
\nopagebreak
\end{proof}

Let now $\KKK$ be the kernel bundle on $\Xi$
whose fiber at $x \in \Xi$ is $\ker x$. Thus $\KKK$ is a
rank 2 vector subbundle of the trivial rank 8 vector bundle
$E_\Xi =E\otimes_\CC \OOO_\Xi$ over $\Xi_0=\Xi\setminus\Omega$. 
Let
$\phi :X\lra H\cap\Xi$ be a  representation of a nonsingular quartic 3-fold
$X\subset\PP^4$ as a linear section of $\Xi$. Giving $\phi$ is
equivalent to specifying a skew-symmetric $8\times 8$ matrix $M$
of linear forms
such that the equation of $X$ is $\Pf (M)=0$. Such a representation
yields a rank 2 vector bundle 
$\EEE =\EEE_\phi$ on $X$, defined by $\EEE =\phi^*\KKK$.
According to \cite{B}, Proposition 8.2, the scheme of zeros of
any section $s\neq 0$ of $\EEE$ is an arithmetically Cohen--
Macaulay  (ACM) 1-dimensional scheme $C$ of degree 14,
not contained in any quadric hypersurface
and such that its canonical bundle $\omega_C\simeq \OOO_C(2)$.
Varieties satisfying the last condition are usually called half-canonical.
Moreover, $\EEE$ is also ACM and has
a resolution of the form
\begin{equation}\label{EACM}
0\lra\OOO_{\PP^4}(-1)^8\mapto{M}\OOO_{\PP^4}^8\lra\EEE\lra 0
\end{equation}
This implies in particular that two Pfaffian representations $\phi_1,
\phi_2$ are equivalent if and only if the corresponding vector
bundles $\EEE_1,\EEE_2$ are isomorphic.  By (8.1) ibid., $\EEE$ can be
given also by Serre's construction as the middle term of the
extension
\begin{equation}\label{serre}
0\lra\OOO_X\lra \EEE\lra \III_{C,X}(3)\lra 0\ ,
\end{equation} 
where $\III_{C,X}$ denotes the ideal sheaf of $C$ in $X$.
Thus, the following assertion holds.

\begin{corollary}\label{7-dim}
A generic quartic $3$-fold $X\subset\PP^4$ has a $7$-dimensional family
of isomorphism classes of rank $2$ ACM vector bundles $\EEE$ with
$\det\EEE\simeq\OOO (3)$ and $h^0(\EEE )=8$
which are characterized by either one of the following equivalent properties:

(i) $\EEE$ as a sheaf on $\PP^4$ possesses a resolution of the form (\ref{EACM})
with a skew-symmetric matrix of linear forms $M$;

(ii) the scheme of zeros of
any section $s\neq 0$ of $\EEE$ is an ACM half-canonical curve $C$ of degree $14$
and arithmetic genus $15$, not contained in any quadric hypersurface in $\PP^4$;

(iii) $\EEE$ can be obtained by Serre's construction from a curve $C\subset X$
as in (ii).
\end{corollary}

In fact the vector bundles under consideration
are stable, so the above 7-parameter family is a part
of the moduli space of vector bundles.

\begin{theorem}\label{M202}
Let $X$ be a generic quartic $3$-fold and
$M_X(2;-1,6)$ the moduli space  of stable rank $2$ vector bundles
$\GGG$ on $X$ with Chern classes $c_1=-[H],c_2=6[l]$, where $[l]\in H^2(X,\ZZ )$
is the class of a line. Then the isomorphism classes of the ACM vector bundles
of the form $\GGG =\EEE (-2)$,
where $\EEE$ are vector bundles introduced in Corollary \ref{7-dim}, form an
irreducible open subset $M_X$ of dimension $7$
in the nonsingular locus of $M_X(2;-1,6)$.
\end{theorem}

\begin{proof} {\em Stability}.  If $\EEE$ is given by the extension
(\ref{serre}), then twisting by $\OOO_X(-2)$ and using
$h^0(\III_{C,X}(k))=0$ for $k\leq 2$ ((ii) of Lemma \ref{7-dim}),
we see that $h^0(\EEE (-2))=0$. The stability follows from Lemma \ref{criterion}.

{\em Smoothness and dimension}. The stability \sloppy implies that $\EEE$ is
simple, that is $h^0(\EEE^\dual\otimes\EEE )=1$. Hence the tangent space
$T_{[\EEE ]} M_X(2;0,2[l])$ at $[\EEE ]$ is identified with
$\Ext^1(\EEE ,\EEE )=H^1(X, \EEE^\dual\otimes\EEE )$, and if
$H^2(X, \EEE^\dual\otimes\EEE )=0$, then $M_X(2;0,2[l])$ is smooth
at $[\EEE ]$ of local dimension $\dim_{[\EEE ]} M_X(2;0,2[l])=
h^1(\EEE^\dual\otimes\EEE )$.

As $\mt{rk}\EEE=2$, we have $\EEE^\dual\simeq\EEE\otimes(\det\EEE )^{-1}\simeq
\EEE (-3)$. By Serre duality, $h^3(\EEE^\dual\otimes\EEE )=
h^0(\EEE^\dual\otimes\EEE (-1))=0$. 
By (\ref{EACM}), 
$h^0(\EEE (-3))=\chi (\EEE (-3))=0$.
Together with the ACM property for
$\EEE$ this gives $h^i(\EEE (-3))=0$ for all $i\in\ZZ$.
Now, from (\ref{serre})
tensored by $\EEE (-3)$, we obtain the isomorphisms
\begin{equation}\label{EIC}
H^i(\EEE^\dual\otimes\EEE )=H^i(\EEE\otimes \EEE (-3))=
H^i(\EEE\otimes\III_C)\ \ \ \forall\ i\in\ZZ\ .
\end{equation}
Further, the restriction sequence
\begin{equation}\label{restriction}
0\lra \EEE\otimes\III_C\lra\EEE\lra\EEE |_C\lra 0
\end{equation}
yields $\chi (\EEE\otimes\III_C)=\chi (\EEE )-\chi (\EEE |_C)=
8-14=-6$, so to finish the proof, it remains to prove the
vanishing of $h^2(\EEE\otimes\III_C)$. By \cite{B}, (8.9), the
vanishing of $h^2(\End_0 (\EEE))$ follows
from the fact that
the map $f$, introduced in the proof
of Proposition \ref{a-la-Adler}, is dominant.
As $h^2(\OOO_X)=0$, we have $0=h^2(\End_0 (\EEE))=h^2(\End (\EEE))=
h^2(\EEE\otimes\III_C)$, and we are done.
\end{proof}

\section{Curves of degree 14 and genus 15 in $\PP^4$}
\label{curves}

Let $X=\{ F=0\}$ be a generic quartic 3-fold in $\PP^4$,
and $X=H\cap\Xi$ (so, the $\PP^4$
is identified with $H$) a Pfaffian representation for $X$. For the sake
of functoriality, we should have defined $\Xi$ as embedded in
$\PP(\wedge^2(E^\dual))$, so that the points $x\in X$ be interpreted
as alternating bilinear forms of rank 6 on $E$, whilst $G=G(2,8)\subset
\PP(\wedge^2E)$; to avoid this dichotomy we will
work in coordinates, identifying $E$ with $E^\dual$. Let $\EEE$ be the
corresponding rank 2 vector bundle and $C$ the scheme of zeros of a section
$s\neq 0$ of $\EEE$. Let $\Hff$, resp. $\Hff^X$ denote the union
of the components of the Hibert scheme of curves in $\PP^4$, resp. $X$ whose
generic points represent a curve $C$ as above. For generic $s$,
the curve $C$ is nonsingular.

Similarly to the previous section, introduce
the rank filtration on the $7\times 7$ skew-symmetric matrices:
$G'=G(2,7)\subset Z\subset\PP^{20}=\PP(\wedge^2(\CC^7))$. According to
\cite{R}, we have $\dim G'=10$, $\deg G'=42$, $\omega_{G'}=\OOO_{G'}(-7)$,
$\dim Z=17$, $\deg Z = 14$, $\omega_Z=\OOO_Z(-14)$. $G'$ will be identified
with a subvariety of $G$ for the standard inclusion $\CC^7\subset\CC^8$.

\begin{proposition}\label{hilb}
The following assertions hold:

(i) $h^0(\NNN_{C/X})=14$, $h^1(\NNN_{C/X})=0$. Hence $\Hff^X$ is smooth
at $[C]$ of local dimension $14$.

(ii) $h^0(\NNN_{C/\PP^4})=56$, $h^1(\NNN_{C/\PP^4})=0$. Hence $\Hff$ is smooth
at $[C]$ of local dimension $56$.

(iii) $C$ can be identified with a section of the rank $4$ locus $Z$
of $7\times 7$ skew-symmetric matrices
by a $4$-dimensional linear subspace $L\subset\PP^{20}$.
\end{proposition}

\begin{proof}
(i) 
The restriction sequence (\ref{restriction})
yields $h^2(\EEE\otimes\III_C)=h^1(\EEE |_C)$. We proved in Theorem \ref{M202}
the vanishing of $h^2(\EEE\otimes\III_C)$.
As $C$ is the
scheme of zeros of a section of $\EEE$, we have $\EEE |_C\simeq
\NNN_{C/X}$. So, we obtain  $h^1(\NNN_{C/X})=0$. By Riemann--Roch,
$h^0(\NNN_{C/X})=14$ and we are done.

(ii) We have $h^1(\NNN_{C/X})=0$. We are going to show that this implies
$h^1(\NNN_{C/\PP^4})=0$. First, by Serre duality,
$0=h^1(\NNN_{C/X})=h^0(\NNN^\dual_{C/X}(2))$. From the restriction sequence
$$
0\lra\III_{C,X}(2)\lra\OOO_X(2)\lra\OOO_C(2)\lra 0
$$
and from the fact that $\omega_C=\OOO_C(2)$, we deduce
that $h^1(\III_{C,X}(2))=0$. Now, the exact triple
$$
0\lra\III^2_{C,X}(2)\lra\III_{C,X}(2)\lra\NNN^\dual_{C/X}(2)\lra 0
$$
yields $h^1(\III^2_{C,X}(2))=0$.
The ACM property for $C$ and 
$$
0\lra\III_{C,\PP^4}(-2)\mapto{F}\III^2_{C,\PP^4}(2)\lra
\III^2_{C,X}(2)\lra 0
$$
imply $h^1(\III^2_{C/\PP^4}(2))=0$. Now, the triple
$$
0\lra\III^2_{C,\PP^4}(2)\lra\III_{C,\PP^4}(2)\lra\NNN^\dual_{C,\PP^4}(2)\lra 0
$$
and the Serre duality give $h^0(\NNN^\dual_{C/\PP^4}(2))=h^1(\NNN_{C/\PP^4})=0$.
By Riemann--Roch, $h^0(\NNN_{C/\PP^4})=56$.

(iii) The sections of $\EEE$ are naturally identified
with elements of $E^\dual$ via the embedding of $\EEE$ into
the trivial rank 8 vector bundle $E_X=E\otimes\OOO_X$.
Let $\Cl :\Xi\setminus\Omega\lra G=G(2,8)$ be the classifying map, sending
each $x\in \Xi\setminus\Omega$ to the projectivized kernel
of $x$, considered as a point
of $G$, and $\Cl_X$ the restriction of $\Cl$ to $X$. 
We can choose the coordinates in $E$ in such a way that
$s=x_7$. Hence $C=\Cl_X^{-1}(\si_{11}(\PP^6))$, where $\PP^6$ is the
hyperplane $\{ x_7=0\}$ in $\PP^7=\PP (E)$, and $\si_{11}(\PP^6)
= G'\subset G$ is the Schubert subvariety of all the lines
contained in the hyperplane. We can also write
$C=\Cl^{-1}(G')\cap H$.
The closure of the $24$-fold ${\Cl}^{-1}(G')$ in $\Xi$
is defined by the $7$ cubic Pfaffians
$\Pf_{r7}(x), 0 \le r \le 6$.

As cubic forms, the Pfaffians $\Pf_{r7}(x)$, $0 \le r \le 6$
do not depend on the variables $x_{p7}$, $0 \le p \le 7$.
Therefore ${\Cl}^{-1}(G')$ is isomorphic
to the cone $C(Z)$ with vertex (or ridge)
$\overline{\PP}^6 = <e_{07},...,e_{67}>$
and base

$Z =
\{ x': \Pf_{07}{x'} =\ldots = \Pf_{67}{x'} = 0 \}
\subset
{\PP}({\wedge}^2 \ <e_0,...,e_6>)$;

\noindent here $x' = (x_{pq})_{0 \le p,q \le 6}$
is the 8-th principal adjoint matrix of the matrix
$x$, i.e. $x'$ is obtained from $x$ by deleting
its last column and row. It is well known that the
vanishing of the principal minors of order $2n$ of a skew-symmetric 
$(2n+1)\times (2n+1)$ matrix
is equivalent to the vanishing of all its minors of order $2n$,
so $Z$ is the locus of $7\times 7$ skew-symmetric matrices
of rank $4$. The projection $\pi :\PP^{27}\dasharrow
\PP^{20}$ with center $\overline{\PP}^6$ maps isomorphically
(for generic $H$) the intersection $H\cap C(Z)$ to
$L\cap Z$, where $L=\pi (H)$. This ends the proof.
\end{proof}

Let $\MMM_{g}$ denote the
moduli space of smooth curves of genus $g$ and $\MMM_g^r$ the subvariety
of $\MMM_g$ parametrizing half-canonical curves with a theta-characteristic
$D$ such that $\dim |D|=r$.

\begin{corollary} \label{M154}
The following assertions hold:

(i) $\Hff$ is irreducible of dimension $56$.

(ii) For generic $\LLL\in \mt{Lin}(\PP^4, \PP^{20})$,
the stabilizer of $\LLL$ in $PGL(7)$, acting on the
right, is finite, and the natural map
$\mt{Lin}(\PP^4, \PP^{20})/PGL(7)\dasharrow\Hff$ is generically finite.

(iii) The natural map $g:{\raisebox{-2pt}{$PGL(5)$}\!\!\setminus\!\!
\raisebox{2pt}{$\mt{Lin}(\PP^4,\PP^{20})$}\! /
\!\raisebox{-2pt}{$PGL(7)$}}\dasharrow
\MMM_{15}^4$ is generically finite and its image is a $32$-dimensional 
irreducible component $\Mfo$ of $\MMM_{15}^4$.
\end{corollary}

\begin{proof}
(i) Indeed, $\Hff$ is the image of $\mt{Lin}(\PP^4, \PP^{20})$.

(ii) This follows from the count
of dimensions: $\dim \mt{Lin}(\PP^4, \PP^{20})-\dim PGL(7)=(5\cdot 21-1)-
(7^2-1)=56=\dim\Hff$.

(iii) According to Harris \cite{H}, if $r\leq\frac{1}{2}(g-1)$,
then the codimension of any component of $\MMM_g^r$ in $\MMM_g$
is at most $\frac{1}{2}r(r+1)$. Applying this to our case, we
see that the dimension of every component of $\MMM_{15}^4$ is
at least 32. Hence the component $\Mfo$, containing the image
of $\Hff$, is of dimension $\geq 32$.
The dimension of ${\raisebox{-2pt}{$PGL(5)$}\!\!\setminus\!\!
\raisebox{2pt}{$\mt{Lin}(\PP^4,\PP^{20})$}\! /
\!\raisebox{-2pt}{$PGL(7)$}}$ is 32, so
it remains to show that $g$ is dominant over $\Mfo$.

Take a generic $C$ from the image of $g$. $C$ is a smooth ACM
curve in $\PP^4$. 
By the definition of $\MMM_g^r$, every small (analytic or \'etale)
deformation of $C$
is accompanied by a deformation of the theta-characteristic $D$
embedding $C$ into $\PP^4$. The ACM property being generic,
any generic small deformation of $C$ is again in the image of $g$,
and we are done.
\end{proof}

\begin{remark}\label{threemat}
In (ii) of the lemma, the stabilizer $G_\LLL$ of $\LLL$
might act by non-trivial automorphisms of $C$. As $\Aut (C)$ is
finite, the subgroup $H_\LLL$ fixing pointwise $C$, and hence $L=\LLL (\PP^4)$,
is of finite index in $G_\LLL$. So, the first assertion of (ii)
is equivalent to saying that $H_\LLL$ is finite.
One can strengthen this assertion: the subgroup of $PGL(2n+1)$ fixing pointwise
a generic linear $\PP^2\subset\PP (\wedge^2\CC^{2n+1})$ for $n\geq 2$ is finite. This
is easily reduced to the $2n$-dimensional case, stated in \cite{B}, (5.3).
\end{remark}

\begin{proposition}\label{P7bundle}
Let $\Hfo^X\subset\Hff^X$ be the locus of ACM half-canonical curves 
$C\subset X$ of degree $14$
and arithmetic genus $15$, not contained in any quadric hypersurface in $\PP^4$,
and $M_X\subset M_X(2;0,2[l])$ the open set defined in Theorem \ref{M202}.
Then the Serre construction  defines a morphism $\phi :\Hfo^X\lra M_X$
with fiber $\PP^7$. Moreover, $\Hfo^X$ is isomorphic locally in the \'etale
topology over $M_X$ to a projectivized rank $8$ vector bundle on $M_X$.
\end{proposition}

\begin{proof}
It is easily seen that $\dim\Ext^1 (\III_C(3),\OOO_X)=1$, so, given $C$, 
the Serre construction determines $\EEE$ uniquely. This yields $\phi$
as a set theoretic map. An obvious relativization of the Serre construction
shows that it is indeed a morphism.

Further, we have $h^0(\EEE\otimes \III_C)=1$
by stability of $\EEE$ and (\ref{EIC}), so the projective space $\PP^7
=\PP (H^0(\EEE ))$ is injected into $\Hff^X$ by sending each section $s\neq 0$
of $\EEE$ to its scheme of zeros. Hence the fibers of $\phi$ are set-theoretically
7-dimensional projective spaces. The proof of the last assertion of the proposition
is completely similar to that of Lemma 5.3 in \cite{MT}.
\end{proof}

\section{Abel--Jacobi map}

We are going to remind briefly the Clemens--Griffiths technique
for the calculation of the differential of the Abel--Jacobi map,
following Welters \cite{W}, Sect. 2. Let $X$ be a nonsingular
projective 3-fold with $h^{03}=0$,
and $X\subset W$ an embedding in a nonsingular possibly non-complete
4-fold. Let $\Phi :B\lra J^2(X)$ be the Abel--Jacobi map, where $B$ is
the base of a certain family of curves on $X$. The differential
$d\Phi_{[Z]}$ at \sloppy a point $[Z]\in B$, representing a curve $Z$,
factors into the composition of
the infinitesimal classifying map $T_{B,b}\lra H^0(Z,\NNN_{Z/X})$
and of the universal ``infinitesimal Abel--Jacobi map"
$\psi_Z:H^0(Z,\NNN_{Z/X})\lra 
H^1(X,\Omega^2_X)^\dual=T_{J_1(X),0}$. The adjoint $\psi_Z^\dual$
is identified by the following commutative square:

\begin{equation}\label{CDWelters}
\begin{CD}
H^0(X,\NNN_{X/W}\otimes\omega_X) @>{R}>> H^1(X,\Omega^2_X) \\
@V{r_Z}VV @VV{\psi_Z^\dual}V \\
H^0(Z,\NNN_{X/W}\otimes\omega_X|_Z) @>{\beta_Z}>> H^0(Z, \NNN_{Z/X})^\dual.\\
\end{CD}
\end{equation}
Here $r_Z$ is the map of restriction to $Z$, and the whole square
(upon natural identifications)
is the $H^0\rar H^1$ part of the commutative diagram of the long exact cohomology
sequences associated to the following commutative diagram
of sheaves:

\begin{equation}\label{CDsheaves}
\begin{array}{ccccccccc}
\scriptstyle{ 0} &\scriptstyle{ \rar} &\scriptstyle{ \Omega^2_X} &\scriptstyle{ \rar}
&        \scriptstyle{ \Omega^3_W\otimes\NNN_{X/W} }&\scriptstyle{ \rar}
&                  \scriptstyle{  \Omega^3_X\otimes\NNN_{X/W} }&\scriptstyle{ \rar }&\scriptstyle{ 0 }\\
& & \downarrow & &\downarrow & & \downarrow & & \\
\scriptstyle{ 0} &\scriptstyle{\rar} &\scriptstyle{ \Omega^3_X\otimes \NNN_{Z/X} }&\scriptstyle{ \rar} &\scriptstyle{ \Omega^3_X\otimes\NNN_{Z/W}}
      &  \scriptstyle{    \rar} &\scriptstyle{ \Omega^3_X\otimes\NNN_{X/W}\otimes\OOO_Z} &\scriptstyle{ \rar} &\scriptstyle{ 0}
%\end{CD}
\end{array}
\end{equation}

Specifying all this to the case when $X$ is a generic quartic 3-fold,
$Z=C\subset X$ a generic curve from $\Hff^X$, $W=\PP^4$, we see
that the dimensions in (\ref{CDWelters}) form the array
$\left(\begin{array}{cc}\scriptstyle 35&\scriptstyle 30\\ \scriptstyle 28&\scriptstyle 14
\end{array}\right)$, that $R, r_C$ are surjective
and that $\corank \beta_C=\corank \psi_C^\dual=h^1(\NNN_{C/\PP^4}(-1))$.
Dualizing, we obtain:
\begin{lemma}\label{h-one}
For $C, X$ as above, \sloppy $\dim\ker\psi_C=h^1(\NNN_{C/\PP^4}(-1))$,
$\dim\im\psi_C=14-h^1(\NNN_{C/\PP^4}(-1))$.
\end{lemma}

\noindent We have $\chi (\NNN_{C/\PP^4}(-1))=14$, hence
$h^0 (\NNN_{C/\PP^4}(-1))=14+h^1(\NNN_{C/\PP^4}(-1))$.

\begin{lemma}\label{h-zero}
$h^0(\NNN_{C/\PP^4}(-1))=21$.
\end{lemma}

\begin{proof}
Twisting the 4 exact triples in the proof of Proposition \ref{hilb} by
$\OOO (1)$, one can see that
the assertion is equivalent to
$$
h^2(\III_{C,\PP^4}^2(3))=21\ ,\ \ h^i(\III_{C,\PP^4}^2(3))=0\ \ \forall\ \ i\neq 2 .
$$
The last equalities follow immediately from the resolution
for $\III_{C,\PP^4}^2(3)$, obtained from
(4) of \cite{R} by restriction to $L=\PP^4\subset\PP^6$ and twisting by $\OOO (3)$: 
$$
0\rar 21\OOO_{\PP^4}(-5)\rar 48\OOO_{\PP^4}(-4)\rar 28\OOO_{\PP^4}(-3)\rar 
\III_{C,\PP^4}^2(3)\rar 0.
$$
\end{proof}

\begin{remark}
One can interprete the elements of $H^0(\NNN_{C/\PP^4}(-1))$
as infinitesimal deformations of $C$ preserving 14 points
of some hyperplane section $S=C\cap h$ of $C$.
It is easy to understand this value geometrically in constructing explicitly a 
21-dimensional family of global, non-infinitesimal
deformations of $C$ which preserve $S$; one can show that every
1-parameter infinitesimal deformation lifts to a global one at list for
generic $C,h$.

Indeed, lift $C$ to an element $A\in \mt{Lin}(\PP^4,\PP^{20})$,
$C=A^{-1}(L\cap Z)$, $L=A(\PP^4)$. The ACM
property of $C$ implies that $S$ spans $h$, so the set $U$ of
global deformations $A'$ of $A$ such that $S\subset C'=
{A'}^{-1}(A'(\PP^4)\cap Z)$ are exactly the elements
$A'$ with the property $A|_h= A'|h$. Identify
$\mt{Lin}(\PP^4,\PP^{20})$ with the open subset of the Grassmannian
$G(5,26)$ parametrizing the graphs of the linear injective maps
from $\CC^5$ to $\CC^{21}$. The graphs of the above elements $A'$
correspond to those 4-dimensional planes in $\PP^{25}$ which
contain a fixed $\PP^3$, the graph of $A|_h$. Thus $U$ is identified with an
open subset of $\PP^{21}\subset G(5,26)$.

We can assume $A,C, h$ generic, so that $A(h)$ is a generic linear $\PP^3$
in $\PP^{20}$. The ACM property for $C$ implies that the 14 points $S$
are in a sufficiently general position, so that the stabilizer of
$S$ in $PGL(4)=\Aut (h)$ is finite, and the subgroup fixing $S$ pointwise is trivial.
This observation and Remark \ref{threemat} imply that
the stabilizer of $A(h)$ in $PGL(7)$ is finite. Hence the orbits of $PGL(7)$
have only finite intersections with $U$. Hence the map of $U$ to the
quotient by $PGL(7)$ is quasi-finite, as well as that to $\Hff$,
and its differential is injective at (a generic) $A$.

We have obtained a 21-dimensional family of global deformations of $C$ preserving $S$.
Now we want to show that any 1-parameter infinitesimal deformation
of $C$ can be lifted to a global 1-parameter one in the image of $U$.
Indeed, the injectivity of the differential allows to lift the infinitesimal
deformation to $U$. An element of $U$ is a proportionality class of
a $5\times 21$ matrix, and $U$ is an open subset in a linear $\PP^{21}$
inside the projective space of the proportionality classes of
$5\times 21$ matrices, so any infinitesimal deformation in $U$ is obviously
lifted to a linear pencil.
\end{remark}

Lemmas \ref{h-one}, \ref{h-zero} imply that the Abel--Jacobi map $\Phi$ has
a \sloppy 7-dimensional image in the 30-dimensional intermediate Jacobian $J^2(X)$
and 7-dimensional fibers. We can easily identify the irreducible
components of the fiber. Indeed, by Proposition \ref{P7bundle},
each $C$ is contained in a $\PP^7=\PP (H^0(\EEE ))\subset\Hff^X$.
Any rationally connected variety is contracted
by the Abel--Jacobi map, so
each one of its fibers is a union of these $\PP^7$'s. As the dimension of
the fiber is 7, the $\PP^7$'s are irreducible components of the fiber.
Being fibers of $\phi$, the irreducible components do not meet each other, so they
are in fact connected components.
Thus we have proved the following
theorem.

\begin{theorem}\label{AJ}
Let $X$ be a generic quartic $3$-fold. 
Let $\Hfo^X\subset\Hff^X$ be defined as in Proposition \ref{P7bundle}, and
$\Phi :\Hfo^X\lra J^2(X)$ the Abel--Jacobi map. Then the dimension of any component of
$\Phi (\Hfo^X)$ is equal to $7$
and the fibers of $\Phi$ are the unions of finitely many disjoint $7$-dimensional
projective spaces. The natural map $\psi :M_X\lra J^2(X)$, defined by
$\Phi =\psi\circ\phi$, is quasi-finite and \'etale on $M_X$.
\end{theorem}

We get immediately the following obvious corollary:

\begin{corollary}
Every component of $M_X$ has non-negative Kodaira dimension.
\end{corollary}

\nopagebreak
\end{document}